# Statistical Evidence Measured on a Properly Calibrated Scale Across Nested and Non-nested Hypothesis Comparisons

Veronica J. Vieland[1,2], Sang-Cheol Seok[1]

[1]    Battelle Center for Mathematical Medicine

The Research Institute at Nationwide Children's Hospital

575 Children's Crossroad, Columbus OH 43215

sang-cheol.seok@nationwidechildrens.org

[2]    Departments of Pediatrics and Statistics

The Ohio State University

Columbus OH 43215

Author to whom correspondence should be addressed:

Veronica Vieland, veronica.vieland@nationwidechildrens.org

## Abstract

Statistical modeling is often used to measure the strength of evidence for or against hypotheses on given data. We have previously proposed an information-dynamic framework in support of a properly calibrated measurement scale for statistical evidence, borrowing some mathematics from thermodynamics, and showing how an evidential analogue of the ideal gas equation of state could be used to measure evidence for a one-sided binomial hypothesis comparison ("coin is fair" versus "coin is biased towards heads"). Here we take three important steps forward in generalizing the framework beyond this simple example. We (1) extend the scope of application to other forms of hypothesis comparison in the binomial setting; (2) show that doing so requires only the original ideal gas equation plus one simple extension, which has the form of the Van der Waals equation; (3) begin to develop the principles required to resolve a key constant, which enables us to calibrate the measurement scale across applications, and which we find to be related to the familiar statistical concept of degrees of freedom. This paper thus moves our information-dynamic theory substantially closer to the goal of producing a practical, properly calibrated measure of statistical evidence for use in general applications.

**Keywords**: Statistical evidence; information dynamics; thermodynamics

**Introduction**

Statistical modeling is used for a variety of purposes throughout the biological and social sciences, including hypothesis testing and parameter estimation among other things. But there is also a distinct purpose to statistical inference, namely, measurement of the strength of evidence for or against hypotheses in view of data. This is arguably the predominant use of statistical modeling from the point of view of most practicing scientists, as manifested by their persistence in interpreting the p-value as if it were a measure of evidence despite multiple lines of argument against such a practice.

   In previous work, we have argued that for any measure of evidence to be reliably used for scientific purposes, it must be properly *calibrated*, so that one "degree" on the measurement scale always refers to the same amount of underlying evidence, within and across applications [1-3]. Towards this end, we proposed adapting some of the mathematics of thermodynamics as the basis for an absolute (context-independent) measurement scale for evidence [4]. The result was a new theory of information-dynamics, in which different types of information are conserved and interconverted under principles that resemble the 1st two laws of thermodynamics, with evidence emerging as a relationship among information types under certain kinds of transformations [5]. As we argued previously, this provides us both with a formal definition of statistical evidence and with an absolute scale for its measurement, much as thermodynamics itself did for Kelvin's temperature T. But unless this new theory can produce something useful, it is purely speculative and not really a theory at all in the scientific sense, so much as an overgrown analogy.

   Until now, though, the theory has been too limited in scope to be of any practical use, for four reasons. (i) We have previously worked out a concrete application only for a simple coin-tossing model, and we speculated that extension to other statistical models (i.e., forms of the likelihood other than the binomial) might require derivation of a new underlying equation of state (EqS, that is, the formula for computing the evidence E; see below for details) for every distinct statistical model. The principles for deriving these new equations remained, however, unclear. (The text will be simplified by the introduction of a number of abbreviations, of which "EqS" is the first. To assist the reader, abbreviations are summarized in Table 1.) (ii) The original EqS also contained two constants, which we speculated might relate to calibrating evidence measurement across different statistical models, but again, the principles under which the constants could be found were unknown, rendering the issue of calibration across applications moot. (iii) Furthermore, the theory appeared to work correctly only in application to a one-sided hypothesis comparison ("coin is fair" versus "coin is biased towards heads"), failing even for the seemingly simple extension to a two-sided comparison ("coin is fair" versus "coin is biased in either direction"). (iv) Because we depended heavily on the arithmetic of thermodynamics in justifying some components of the theory, it was unclear how to move to general applications without relying upon additional equations to be borrowed from physics and applied in an *ad hoc* manner to the statistical problem. While this issue was ameliorated by the introduction of wholly information-based versions of the 1st and 2nd laws of thermodynamics [5], it remained a concern, particularly in view of our inability to extend the theory beyond the one-sided binomial application.

**Table 1** Summary of Abbreviations

| Abbreviation | Full Name | Description |
|---|---|---|
| BBP | Basic Behavior Pattern | a characteristic of what we mean by "statistical evidence" that any *measure* of evidence must recapitulate |
| d.f. | Degrees of Freedom | one of two constants in the EqS, used to calibrate E across different forms of HC |
| e | evidence | evidence measured on an empirical (uncalibrated) scale |
| E | Evidence | evidence measured on an absolute (context-independent) scale |
| EqS | Equation of State | used to calculate the evidence from features of the likelihood ratio graph |
| HC | Hypothesis Contrast | the forms of the hypotheses in the numerator and denominator of the LR (nested or non-nested; composite or simple) |
| LR | Likelihood Ratio | P(data | Hypothesis 1)/P(data | Hypothesis 2) |
| S | Evidential Entropy | a particular form of Kullback-Leibler divergence, equal to the max log LR |
| TrP | Transition Point | values of $x/n$ at which evidence switches from supporting one hypothesis to supporting the other |
| V | Volume | area (or more generally, volume) under the LR graph |

Thus we were faced with the question of whether the striking connection we had found between the mathematical description of the dynamics of ideal gases and the mathematics of our simple statistical system was really telling us something useful on the statistical side, or whether, by contrast, we had simply stumbled upon a kind of underlying one-sided binomial representation of the ideal gas model in physics — a model of use neither to physicists nor to statisticians. With the results presented in this paper, however, we take an important step forward in laying this concern to rest. We show below how the theory is readily generalized to support a wider range of statistical applications than had previously been considered, and we make strides in laying out the principles under which both the equations of state and the constants can be resolved. In the process, we continue to see connections to the equations of thermodynamics.

Specifically, we generalize the original theory to address the four limitations mentioned above, albeit still in the context of binomial models. We find that equations of state are governed by the different possible forms of hypothesis contrast (HC). We are then able to extend the original results to other HCs, including the two-sided HC that thwarted our earlier attempts at generalization, by introducing a simple extension of the original EqS. We also show a connection between one of the constants and something closely related to the familiar statistical concept of degrees of freedom.

The remainder of the paper is organized as follows. We first (1) briefly review the key methodological principles and results from earlier work, and we illustrate the problem that arises when we move from one-sided to two-sided hypothesis comparisons. In (2) we group binomial HCs into two major Classes, non-nested and nested, and we show that the HCs in Class I can be handled via the original EqS, while a simple modification of this EqS suffices to handle the Class II HCs. In (3) we consider resolution of a key constant across different hypothesis contrasts, and find that it is related to the statistical concept of degrees of freedom. In (4) we illustrate aspects of the behavior of the resulting evidence measure E within and across HCs.

## 1. Review of previous results and the problem with two-sided hypothesis comparisons

We begin with a high level definition of evidence as a relationship between data and hypotheses in the context of a statistical model. We then pose a measurement question: How do we ensure a meaningfully calibrated mapping between (i) the object of measurement, i.e., the evidence or evidence strength, and (ii) the measurement value? Here the object of measurement cannot be directly observed, but must be inferred based on application of a law or principle that maps observable (computable) features of the data onto the evidence. This is known as a nomic measurement problem [6]. There are precedents for solutions to nomic measurement problems, particularly in physics; measurement of temperature is an example [6].

Our guiding methodological principle is that any measure of evidence must verifiably behave like the evidence, in situations in which such verification is possible. In order to establish basic behavior patterns (BBPs) expected of any evidence measure, we consider a simple model and a series of thought experiments, or appeals to intuition. This enables us to articulate basic operational characteristics of what we mean by "statistical evidence." We then check any proposed measure of evidence to be sure that it exhibits the correct BBPs. As the theory is developed, we are also able to observe new patterns of behavior. These are considered iteratively to assess their reasonableness.

These BBPs play a role here that is similar to the role played in some other treatments of evidence by axioms [7] or "conditions" [8]. However, in our methodology the BBPs themselves only support a measure of evidence *e* on an empirical, rather than absolute scale. Any proper empirical measure *e* must exhibit the BBPs. But as long as the *only* criterion is that *e* exhibits the BBPs, the units of *e* remain arbitrary and they are not necessarily comparable across applications. Thus the BBPs constitute a set of necessary, but not sufficient, conditions for a proper evidence measure.

The primary set of thought experiments used to establish the current set of BBPs involves coin-tossing examples, for which our intuitions are clear and consensus is easy to achieve on key points. (Royall [9] also uses a simple binomial set-up as a canonical system for eliciting intuitions about evidence. However, his use of the binomial is quite different from ours. He appeals to intuition in order to calibrate strength of evidence across applications; we appeal to intuition to establish certain properties we expect evidence to exhibit. In our methodology, calibration is a separate process.) Consider a series of $n$ independent coin tosses of which $x$ land heads and $n$-$x$ land tails. Let the probability that the coin lands heads be $\theta$. And consider the two hypotheses $H_1$: "coin is biased towards tails" ($\theta < \frac{1}{2}$), versus $H_2$: "coin is fair" ($\theta = \frac{1}{2}$). We articulate four BBPs up front. (We have described the thought experiments used to motivate the BBPs in detail elsewhere; see, e.g., [5]. Here we simply summarize the BBPs themselves.)

(i) *Change in evidence as a function of n for fixed x/n* For any fixed value of $x/n$, the evidence increases as $n$ increases. The evidence may favor $H_1$ or $H_2$, depending on $x/n$, but in either case, it increases with increasing $n$. BBP(i) is illustrated in Figure 1(a).

(ii) *Change in evidence as a function of x/n for fixed n* For any fixed $n$, as $x/n$ increases from 0 to ½ the evidence in favor of $H_1$ decreases up to some value of $x/n$, after which it increases in favor of $H_2$. We refer to the value of $x/n$ at which the evidence switches from favoring $H_1$ to favoring $H_2$ as the transition point (TrP). We also expect the value of $x/n$ at which a TrP occurs to shift as a function of $n$, as increasingly smaller departures from $x/n = 0.5$ support $H_1$ over $H_2$. BBP(ii) is also illustrated in Figure 1(a).

(iii) *Change in x/n and n for fixed evidence* In order to maintain constant evidence, as $x/n$ increases from 0 to the TrP, $n$ increases; as $x/n$ continues to increase from the TrP to ½, $n$ decreases. These patterns follow from BBP(i) and BBP(ii). BBP(iii) is illustrated in Fig 1(b).

(iv) *Rate of increase of evidence as a function of n for fixed x/n* The same quantity of new data ($n$, $x$) has a *smaller* impact on the evidence the *larger* is the starting value of $n$, or equivalently, the stronger the evidence is before consideration of the new data. E.g., 5 tosses all of which land tails increase the evidence for $H_1$ by a greater amount if they are preceded by 2 tails in a row, compared to if they are preceded by 100 tails in a row. BBP(iv) is illustrated in Fig 1(c).

We can summarize by saying that the three quantities $n$, $x$, and evidence $e$, enter into an EqS, in which holding any one of the three constant while allowing a second to change necessitates a compensatory change in the third. Here $e$ itself is simply *defined* as the third fundamental entity in the set. At this point no particular measurement scale is assigned to $e$, and therefore numerical values are not assigned to $e$ and $e$-axes are not labeled in the figures. Figure 1 is intended to illustrate behavior *patterns* only, rather than specific numerical results. (In Figure 1 and subsequent Figures, $n$ and $x$ are treated as continuous rather than integer, in order to smooth the graphs particularly for small $n$.)

Our overarching methodological principle is that any proposed measure of evidence must exhibit these basic patterns of behavior. While we are free to use any methods we like to discover or invent a statistical EqS, applying this principle to our simple set of BBPs severely restricts the set of permissible EqSs. For example, it is easily shown that the p-value and –log[p-value], the maximum likelihood ratio and its logarithm, and the Bayes factor all violate one or more of the BBPs; e.g., they all violate BBP(iv). Thus the equations used to calculate these quantities cannot serve as EqSs for measurement of evidence. Of course, in enumerating the BBPs thus far we have considered only single-parameter cases. Generalizations to multi-parameter settings may entail additional considerations.

**Figure 1** Basic Behavior Patterns for evidence *e*: (a) *e* as a function of *x/n* for different values of *n*, illustrating BBPs(i) and (ii) (dots mark the TrP, or minimum point, on each curve); (b) iso-*e* contours for different values of *e*, (higher contours represent larger values of *e*), illustrating BBP(iii); (c) *e* as a function of *n* for any fixed *x/n*, illustrating BBP(iv). Because *e* is on an empirical (relative) measurement scale, numerical values are not assigned to *e* and *e*-axes are not labeled in the figures.

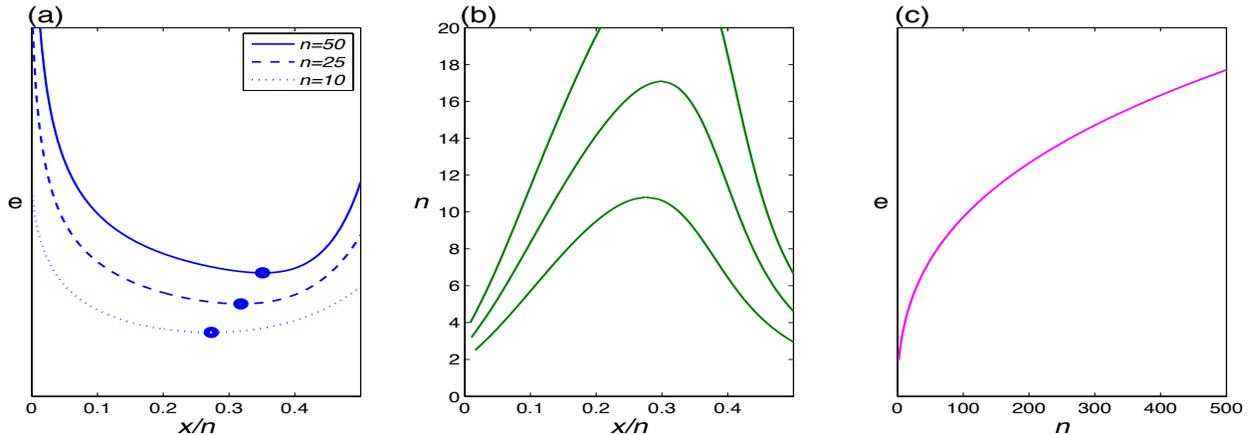

We treat the likelihood ratio (LR) as fundamental. Originally [4] we considered only the special binomial form of LR, $LR(\theta, \theta = 0.5; n, x) = \frac{\theta^x (1-\theta)^{n-x}}{(0.5)^n}$, with $0 \le \theta \le 0.5$ (a composite vs. simple HC with the simple hypothesis on the boundary of the parameter space). The EqS for this set-up was originally derived via the information-dynamic analogue of thermodynamic systems [4]. Here we focus only on the EqS itself and not its derivation. This EqS turned out to be a function of two aspects of the LR (and two constants; see below): (i) the logarithm of the maximum LR, which we treated as an entropy term (see Appendix 1) and denoted as S; and (ii) the area under the LR graph, denoted V, which is related to, though distinct from, the Bayes factor [10] and the Bayes ratio in statistical genetics [11]. Originally these equations were given as

$$S = \log \left[ \frac{\left(\frac{x}{n}\right)^x \left(1-\frac{x}{n}\right)^{n-x}}{\left(\frac{1}{2}\right)^n} \right] \tag{1.1}$$

and

$$V = \int_0^{0.5} LR(\theta, 0.5; n, x) \, d\theta. \tag{1.2}$$

In [4] we derived a simple EqS in the form

$$S = c_1 \log E + c_2 \log V \tag{1.3}$$

for $c_1, c_2$ constants, where E represents evidence measured on an absolute, and not merely an empirical, scale [4]. Equation (1.3) is identical in form to the ideal gas EqS in physics, although we assign different (non-physical) meanings to each of the constituent terms.

Because the focus of this paper is on application of the theory, we do not address the meaning of E in any detail here. But in brief, E is defined as the proportionality between (i) the change in a certain form of information with the influx of new data, and (ii) the entropy, such that the degree of E retains constant meaning across the measurement scale and, given the correct EqS, across applications. See [5] for details.

From (1.3) we have a simple calculation formula for E as

$$E = \left( \frac{exp^S}{V^{c_2}} \right)^{1/c_1}. \tag{1.4}$$

It is readily verified that using (1.4) yields an evidence measure E that exhibits the BBPs described above; in fact, Figure 1 was drawn by applying this equation. In previous work we noted that the principles for determining the constants remained to be discovered, and we set the values somewhat arbitrarily to $c_2 = 1$ and $c_1 = 1.5$. We have found that $c_2 = 1$ is required to maintain the BBPs. We continue to use this value throughout the remainder of this paper, but we have retained $c_2$ in the equations as a reminder that it may become important in future extensions of the theory. We return to resolution of $c_1$ in §3 below.

Figure 2 illustrates the problem that we faced in attempting to use (1.4) for a two-sided hypothesis comparison. For given *n* and viewed as functions of *x/n*, Figure 1(a) and Figure 2(a) exhibit similar shapes. In Figure 1(a) (one-sided comparison) the minimum value of E corresponds to the TrP, the *x*/*n* value at which the evidence begins (reading left to right) to favor $\theta_2 = \frac{1}{2}$. Figure 2(a) might at first appear to be a simple extrapolation, but in fact it must be fundamentally wrong. The minimum value should occur at the TrP, the *x*/*n* value at which the evidence begins to favor $\theta_2 = \frac{1}{2}$. But here the minimum point is occurring at the value *x*/*n* = $\frac{1}{2}$, regardless of *n*. Thus this minimum point no longer has the interpretation of being a TrP, that is, a point at which the evidence starts to favor $\theta_2 = \frac{1}{2}$. Indeed, there is no such thing as evidence in favor of $\theta_2 = \frac{1}{2}$ here, since even as *n* increases the evidence remains at a minimum when the data fit perfectly with H$_2$. Figure 2(b) illustrates the pattern (although not necessarily the actual numbers) we should obtain, which requires two TrPs, one on each side of $\theta_2$. In contrast to Figure 2(a), Figure 2(b) represents the correct reflection of the behavior in Figure 1(a) onto the region *x*/*n* > 0.5. In the following section we show how to adjust the EqS to produce the correct pattern as shown in Figure 2(b).

**Figure 2** The problem with using the original equation of state in application to the two-sided hypothesis contrast: (a) E as a function of *x*/*n* for different values of *n*, using the original EqS, illustrating the absence of a true TrP (dots indicate minimum value of E); (b) the expected pattern of behavior of behavior of *e* in the two-sided case, illustrating the correct TrP behavior, with symmetric TrPs on either side of 0.5, converging towards 0.5 as *n* increases. In (a), because we are using the EqS (1.4) to calculate the evidence, we label the y-axis E; however, because this is the *wrong* EqS here, numerical values of E are not labeled.

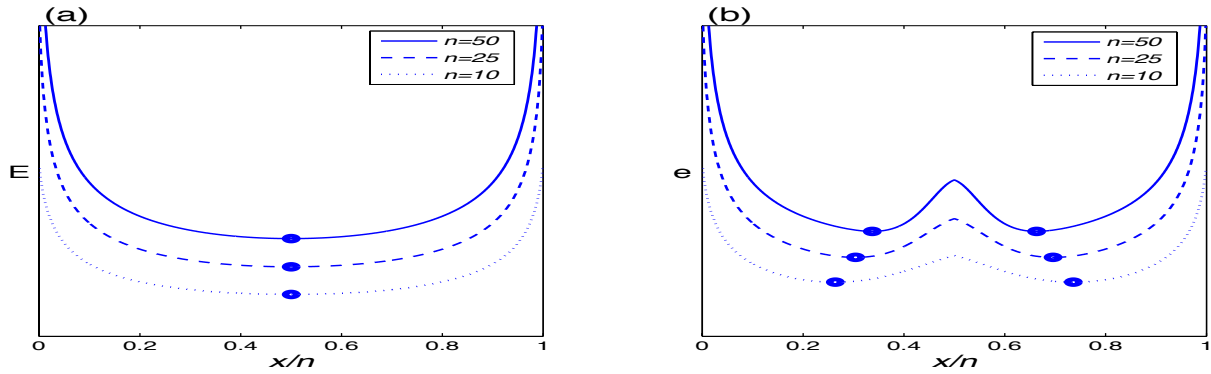

## 2. Equations of state for non-nested and nested HCs

We continue to consider the binomial model with the single parameter θ, and pairs of hypotheses specifying various ranges for θ. We restrict attention to HCs in the form $H_1: \theta \in \Theta_1$ versus $H_2: \theta \in \Theta_2$, where $\Theta_1 \cup \Theta_2 = \Omega$, the set of all possible values of θ. For simplicity of notation, we use subscripts (1, 2) to designate the set of values of θ as stipulated under $H_1$, $H_2$ respectively.

Following a familiar statistical convention, we distinguish two major classes of HC, non-nested (Class I) and nested (Class II). Within each of these classes we can further distinguish (a) composite versus simple HCs from (b) composite versus composite HCs. (Note that our requirement $\Theta_1 \cup \Theta_2 = \Omega$ precludes simple versus simple HCs.) Figure 3 summarizes and illustrates the four resulting HC types.

As shown in the figure, we further restrict attention to HCs in which θ = ½ plays a special role. Specifically, for Class I(a) (the original model [4]), we consider only $H_1: \theta \in [0, \frac{1}{2})$ and $H_2: \theta = \frac{1}{2}$; for Class I(b), we consider only the case $H_1: \theta \in [0, \frac{1}{2}]$ and $H_2: \theta \in (\frac{1}{2}, 1]$; for Class II(a) we consider only the case $H_1: \theta \in [0, 1]$ and $H_2: \theta = \frac{1}{2}$; and for Class II(b) we consider only ranges $\theta_2 \in [\theta_{2l}, \theta_{2r}]$ (where the subscript "*l*" stands for "left" and "*r*" for "right") that are symmetric around the value ½. We have speculated from the start that the unconstrained maximum entropy state of a statistical system, which in the binomial case occurs when θ = ½, plays a special role in this theory. Indeed, in order to maintain the BBPs, calculations have shown that binomial HCs that are not "focused" in some sense on θ = ½ will require further corrections to the underlying EqS. We had also speculated previously that HCs in the form "A vs. not-A" play a special role. Here we extend the theory to include nested hypotheses.

A little thought will show that Class I(b), like the original Class I(a), should have 1 TrP; while Class II(b), like Class II(a), should have 2. The absence of the second TrP was the major reason for feeling that our original EqS did not cover the two-sided case Class II(a), and it appears that Class II(b) will present a similar challenge. Thus we expect both Class I HCs to exhibit the pattern illustrated in Figure 2(a); and both Class II HCs to exhibit the pattern illustrated in Figure 2(b).

Before proceeding we need to generalize our original notation to allow for the additional forms of HC. Let $\hat{\theta} = x/n$, the value of θ that maximizes the likelihood L(θ). Let $\hat{\theta}_i =$ the value of $\theta_i$ (*i* = 1, 2) that maximizes L(θ) within the range imposed by $H_i$.

**Figure 3** Summary and illustration of four basic HCs considered in the text

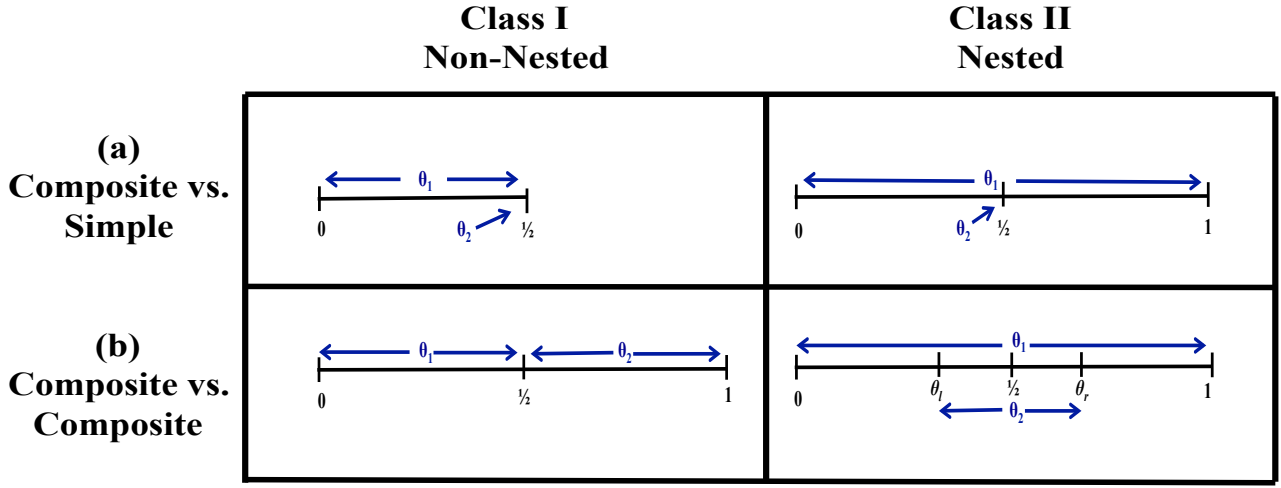

As noted previously, from the start we have viewed S, originally defined as the maximum log LR, as an entropy term; that is, in the original formalism [4] the maximum log LR occupied the place of the term for thermodynamic entropy in the ideal gas EqS. We now explicitly express S as a form of Kullback-Leibler divergence [12] (see Appendix 1 for details). The generalized definition of S (cf. equation (1.1)) becomes

$$S = \sum_{x=0}^{n} P_n(x;\ \hat{\theta}) \log \frac{L(\hat{\theta};\ n,x)}{L(\hat{\theta}_i;\ n,x)}.$$ (2.1)

In the denominator of the LR $i$ =2, except for Class I(b), for which $i = 2$ when $x/n \in \Theta_1$ (i.e., when $x/n \leq \frac{1}{2}$) and $i = 1$ when $x/n \in \Theta_2$ ($x/n > \frac{1}{2}$). In either case $\hat{\theta}_i = \frac{1}{2}$.

We similarly generalize the definition of V (cf. equation (1.2)) to

$$V = \int \frac{L(\theta;\ n,x)}{L(\hat{\theta}_i;\ n,x)}\ d\theta,$$ (2.2)

where for the original one-sided HC, Class I(a), the integral is taken over [0, ½], and for the remaining HCs, the integral is taken over [0,1]. For a simple $\theta_i$, $\hat{\theta}_i = \theta_i$, therefore in application to the original one-sided HC, (2.1) and (2.2) maintain the original definitions for S and V as given in (1.1) and (1.2). From here on, we utilize the generalized definitions of S and V in (2.1) and (2.2).

The original EqS (1.4), which generates all of the correct behaviors for Class I(a), also generates all of the correct behaviors for Class I(b) (see §4 below). Moreover, applying this EqS to Class II(b) *fails* in exactly the same way it does for Class II(a), that is, it fails to generate the second TrP. It turns out that a simple adjustment to (1.4) generates the second TrP for the (symmetric) two-sided binomial HC Class II(a), "coin is fair" versus "coin is biased in either direction." In particular, we adjust our basic EqS by subtracting a term $b$ from V. This yields the new EqS

$$E = \left( \frac{exp^S}{(V-b)^{c_2}} \right)^{1/c_1}. \tag{2.3}$$

The formula for calculating $b$ is given in Appendix 2. Note that while the original EqS (1.4) was in the form of the thermodynamic equation for an ideal gas, (2.3) is in the form of the Van der Waals equation [13].

It is readily verified that (2.3) returns the correct behavior, with two TrPs, as illustrated in Figure 2(b) (indeed, Figure 2(b) was drawn using (2.3)), as well as exhibiting all other BBPs. Equation (2.3) also generates the correct behavior for Class II(b) (see §4). We might have guessed from the outset that the EqS for Class II(b) should be the same as the EqS for Class II(a), since, as the width of the interval $[\theta_{2r} - \theta_{2l}]$ narrows to $\frac{1}{2} \pm \varepsilon$, the two HCs become (approximately) the same, namely, $\theta \neq \frac{1}{2}$ vs. $\theta = \frac{1}{2}$. Therefore for any given data $(n, x)$, as $\varepsilon$ shrinks to 0, they must yield the same value of E. This strongly suggests that a single EqS should govern both types of HC, as indeed turns out to be the case.

### 3. The constant $c_1$ and degrees of freedom

The central point of developing a properly calibrated evidence scale is to be able to meaningfully compare values of the evidence across applications. One obvious way in which we might need to adjust E across different HCs would be to allow for differences in "degrees of freedom" (d.f.). Here we are using d.f. in a generic sense, to signify the difference (or sum, see below) between the dimensionalities of the parameter spaces under the two hypotheses [14]. It should not be surprising if some concept of d.f. enters into the calibration process for E given the familiar role of d.f. in frequentist statistical settings. At the same time, it should also not be surprising if the concept of d.f. enters into our equations *differently* than it does in other statistical settings, due to fundamental differences between our framework and frequentist methodology.

For example, d.f. play a key role in the frequentists' generalized LR $X^2_{d.f.}$ (nested) test. Under broad regularity conditions, familiar mathematics leads to setting the d.f. equal to the difference in the number of parameters being maximized over in the numerator and denominator, respectively, of the maximum log LR. The frequentist d.f. adjustment is required specifically to reflect the fact that the sampling distribution of the maximum LR under the null hypothesis shifts upwards the greater the d.f., and it serves to align Type I error behavior across hypothesis comparisons involving different amounts of maximization. But in our methodology the sampling distribution of the LR is irrelevant. Indeed, we have previously pledged allegiance to a version of the likelihood principle, which is ordinarily understood to preclude consideration of sampling distributions – the distributions of data that might have been but were in fact not observed – when evaluating evidence. (Moreover, the $X^2_{d.f.}$ itself represents *asymptotic* behavior of the maximum log LR. But we are expressly concerned with calibrating evidence measurement in finite samples.) Furthermore, d.f. as a parameter of this particular distribution are applicable only to nested HCs, whereas for us the objective of calibration across applications requires a concept of d.f. that allows a unified treatment of nested and non-nested HCs.

Adherents of the likelihood principle, however, generally eschew any kind of d.f. "correction" to LRs as indicators of evidence strength even in the context of composite hypotheses (see, e.g., [7]). But the premise that a given value of the maximum LR corresponds to the same amount of evidence regardless of the amount of maximization being done strains credulity. Among other problems, this begs the question of overfitting, in which a bigger maximum LR can almost always be achieved by maximizing over additional parameters (up to a model involving one independent parameter for each data point), even in circumstances in which the estimated model can be shown to be getting further from the true model as the maximum LR increases. (See, e.g., the discussion of model fitting versus predictive accuracy in [15]. See also [16], for a coherent pure-likelihoodist resolution of this problem, which avoids "corrections" to the LR for d.f., but which also precludes the possibility of meaningful comparisons of evidence strength across distinct HCs or distinct forms of the likelihood.)

Prior to the new results in §2 above, we had been unable to derive the EqS for HCs other than the original one-sided binomial (with $\theta_2$ on a boundary), and therefore the idea of adjusting the calculation of E to reflect differences in d.f. *across* HCs was moot. But the discovery that just two basic EqSs cover a wide range of HCs strongly suggests that any d.f. adjustment should be captured by some feature of the EqS as shown in (1.4) and (2.3). And, as these equations show, $c_1$ adjusts the magnitude of E for given S and V, which is on the face of it just what we need to do.

It may seem odd to call $c_1$ a constant and then to vary it. We note, however, that in the thermodynamic analogues of our equations (1.4) and (2.3), the position of our $c_1$ is occupied by the physical constant $c_V$, the thermal capacity of a gas at constant volume. This constant varies, e.g., between monatomic and diatomic gases, reflecting the fact that a fixed influx of heat will raise the temperatures of the two gas types by different amounts. Similarly, we can view $c_1$ as a factor that recognizes that different HCs will convert the same amount of new information (or data) into different changes in E. This viewpoint is consonant with our underlying information-dynamic theory [4, 5], which treats transformations of LR graphs in terms of Q (a kind of evidential information influx) and W (information "wasted" during the transformation, in the sense that it does not get converted into a change in E); the sense in which E maintains constant meaning across applications relates specifically to aspects of these transformations (see [4, 5] for details).

The only remaining task then is to find the correct values of $c_1$ for different HCs, as we describe in the following paragraph. Final validation of any specific numerical choices we make at this point regarding $c_1$ will require returning to the original information-dynamic formalism. But we point out here that the choices we have made are far from *ad hoc*. The form of the EqS itself combined with constraints imposed by the BBPs place severe limitations on how values can be assigned to $c_1$ while maintaining reasonable behavior for E within and across HCs.

We have found that we must have $c_1 > 0.5$ in order to maintain BBP(iv). Thus we begin, somewhat arbitrarily but in order to start from an integer value, from a baseline $c_1 = 1.0$. In the case of nested hypotheses ($\Theta_2 \subset \Theta_1$), we add to this baseline value the *sum* $[\theta_{1r} - \theta_{1l}] + [\theta_{2r} - \theta_{2l}]$ of the lengths of the intervals. Heuristically, we sum these lengths because it is possible for $x/n$ to be in either or both intervals simultaneously; thus speaking very loosely, $c_1$ captures a kind of conjunction of the two intervals. In the case of non-nested hypotheses, for which $x/n$ can be in $\Theta_1$ or $\Theta_2$ but not both, a disjunction of the intervals, we add to the baseline the *difference* $[\theta_{1r} - \theta_{1l}] - [\theta_{2r} - \theta_{2l}]$. Using these rules we arrive at $c_1$ = d.f. = 1.5, 1.0, 2.0 and $2 + [\theta_{2r} - \theta_{2l}]$, for Classes I(a), I(b), II(a) and II(b), respectively. Thus for Class II(b), $2 \le c_1 \le 3$. Table 2 shows the assigned values in the context of the EqS for each HC.

**Table 2** Final EqS for each of the four HCs.

| | Class I Non-Nested | Class II Nested |
|---|---|---|
| **(a) Composite vs. Simple** | $E = \left(\dfrac{exp^S}{V}\right)^{1/1.5}$ | $E = \left(\dfrac{exp^S}{V-b}\right)^{1/2}$ |
| **(b) Composite vs. Composite** | $E = \dfrac{exp^S}{V}$ | $E = \left(\dfrac{exp^S}{V-b}\right)^{1/(2+[\theta_{2r}-\theta_{2l}])}$ |

Note that as $c_1$ increases, for given data, E decreases. Thus these values ensure some intuitively reasonable behavior in terms of the conventional role of d.f. adjustments. For instance, for given $x/n$, the two-sided Class II(a) HC will have lower E compared to the one-sided Class I(a) HC, which conforms to the frequentist pattern for one-sided versus two-sided comparisons. We consider the behavior of E in greater detail within and across HCs in §4 below.

**4. Behavior of E within and across HC classes**

Our overarching goal here is to quantify statistical evidence on a common, underlying scale across all four HCs. As noted above, (1.4) and (2.3) ensure the BBPs in application to each HC considered on its own, provided that we set $c_2 = 1$ and $c_1$ as shown in Table 2. In this section we highlight important additional characteristics of E beyond the original BBPs. Some of these characteristics conform to intuitions we had formed in advance, but others constitute newly discovered properties of E – behaviors we did not anticipate, but which nevertheless seem to us to make sense once we observe them.

We begin with Class II(b) on its own, as a function of the size $[\theta_{2r} - \theta_{2l}]$ of the $\Theta_2$ interval. Figure 4 illustrates the behavior of E for Class II(b). Several features of Figure 4 are worth noting. Intuition tells us that for $x/n \approx 0$ or $x/n \approx 1$, as $[\theta_{2r} - \theta_{2l}]$ *increases*, the strength of the evidence in favor of $\theta_1$ should *decrease*, to reflect the fact that even such extreme data represent a smaller and smaller deviation from compatibility with $\theta_2$. This pattern is seen in Figure 4(a), where E = 6.2, 5.4, 4.6 for $[\theta_{2r} - \theta_{2l}]$ = 0.02, 0.20 and 0.40, respectively. For any given $x/n \in \Theta_2$, it also seems reasonable that the evidence, now in favor of $\theta_2$, should *decrease* as $[\theta_{2r} - \theta_{2l}]$ *increases*, again as seen in Figure 4(a). This reflects the fact that $\Theta_2 \subset \Theta_1$, so that evidence to differentiate the two hypotheses is smaller the more they overlap. At the same time, within this interval we would expect $x/n \approx \frac{1}{2}$ to yield the strongest evidence; however, given the overlap between $\Theta_1$ and $\Theta_2$, we would not necessarily expect the evidence at $x/n \approx \frac{1}{2}$ to be substantially larger than the evidence at $x/n$ closer to the $\Theta_2$ boundary. Figure 4(b) illustrates this pattern for different values of n. Note that E is actually maximized at $x/n = \frac{1}{2}$: e.g., for n = 50, at $x/n = \theta_{2l}$ = 0.4, E = 2.75, while at $x/n$ = 0.5, E = 2.78. It is also interesting to note that the TrPs move outward as $[\theta_{2r} - \theta_{2l}]$ increases, as might be expected (Figure 4(a)); while for each fixed $[\theta_{2r} - \theta_{2l}]$, the TrPs are moving inward as *n* increases (Figure 4 (b)), in all cases, converging towards the corresponding left (or right) boundary value of $\theta_2$. Thus in all regards, the adjustment of $c_1$ combined with the Class II EqS seems to yield sensible behavior for E for Class II(b).

**Figure 4** Behavior of E for Class II(b): (a) E as a function of $x/n$ (n = 50) for different ranges for $\theta_2$; (b) E as a function of $x/n$ for $0.4 \leq \theta_2 \leq 0.6$ for different *n*. Note that this graph utilizes the correct EqS. Therefore the y-axis is now labeled as E and numerical values are shown.

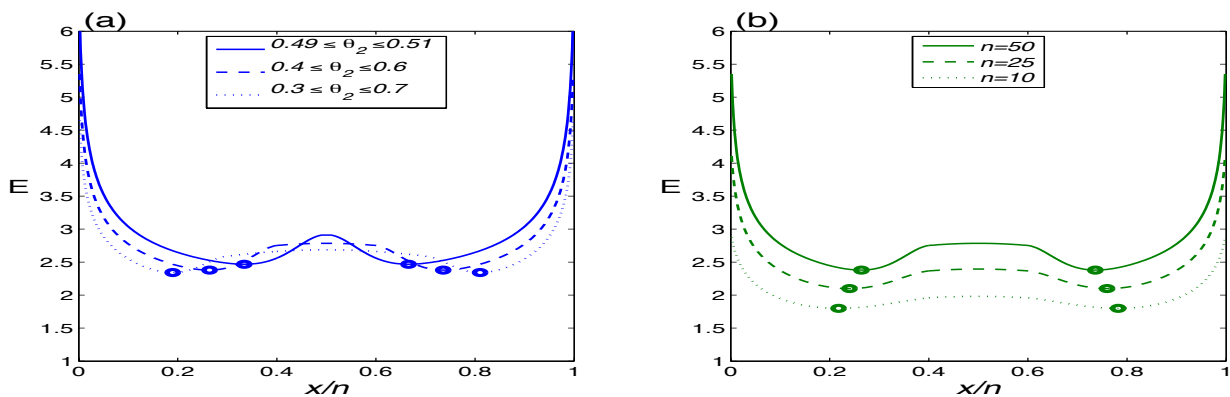

We can also assess the reasonableness of E for Class II(b) in comparison with Class II(a). As $[\theta_{2r} - \theta_{2l}] \rightarrow 0$, $c_1$ becomes the same for Class II(b) and Class II(a), by design. Thus the line in Figure 4(a) representing $0.49 \leq \theta_2 \leq 0.51$ ($c_1$ = 2.02) is virtually identical to what we would obtain under Class II(a) ($c_1 = 2.00$), and for the moment we treat it as a graph of Class II(a). We noted above that that for $x/n \approx 0$ or $x/n \approx 1$, evidence decreases as $[\theta_{2r} - \theta_{2l}]$ increases. We can now see from Figure 4(a) that this also means that evidence is decreasing relative to what would be obtained under a Class II(a) HC. Since under Class II(a) the HC always involves a comparison against $\theta = \frac{1}{2}$, it is reasonable that larger (nested) $[\theta_{2r} - \theta_{2l}]$ would return smaller evidence at these $x/n$ values relative to a comparison against the single value $\theta = \frac{1}{2}$. For $x/n = \frac{1}{2}$, we might have guessed that E in favor of $\theta_2$ should be also smaller for Class II(b) than for Class II(a), since the data are perfectly consistent with both $\theta_1$ and $\theta_2$ but Class II(a) has the more specific $H_2$.

Turning to comparisons across all four HCs, Figure 5 illustrates some additional important behaviors. Across the board, for given $x/n$, E is higher for the Class I HCs than it is for the Class

**Figure 5** Comparative behavior E as a function of $x/n$ ($n = 50$) across all four HCs. For purposes of illustration, $0.4 \leq \theta_2 \leq 0.6$ for Class II(b). TrPs are marked with circles (Class I(a), Class II(a)) or diamonds (Class I(b), Class II(b)).

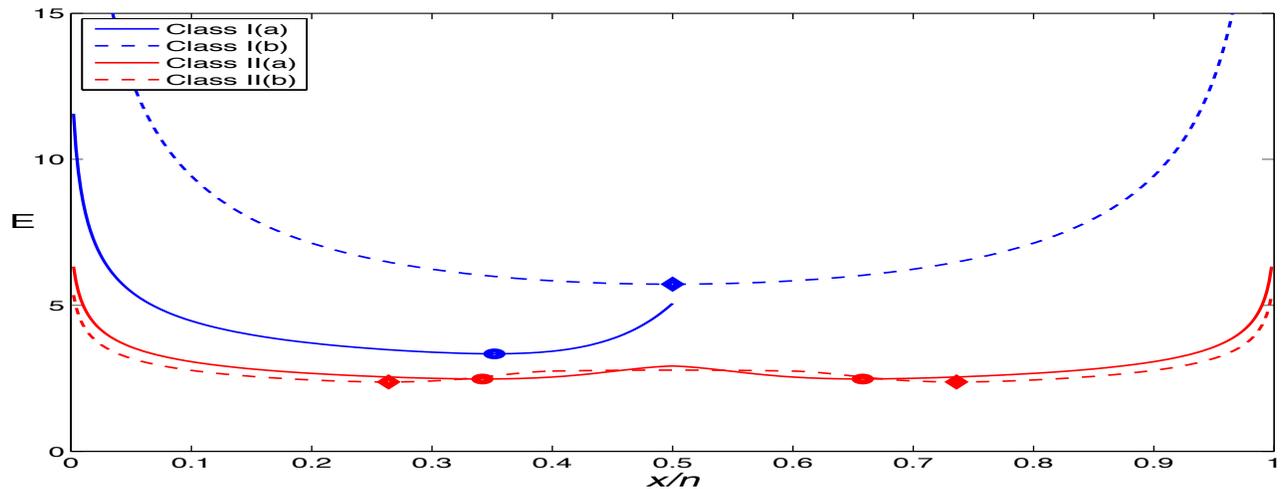

II HCs. This is the result of our assignments for $c_1$, as discussed above, and it makes sense that nested hypotheses would be harder to distinguish compared to non-nested hypotheses for given $n$. Figure 5 also illustrates the relative placement of the TrPs across HCs, which is consistent with, and a generalization of, the BBPs involving the TrP considered in §1 in the context of a single HC. For instance, the TrPs for Class II(b) are further apart than for Class II(a), a pattern we might have anticipated.

Figure 6 reorganizes the representation shown in Figure 5 in terms of "iso-E" contours through the ($n$, $x/n$) space for the different HCs; that is, these graphs display the sets of ($n$, $x/n$) pairs corresponding to the same evidence E. For simplicity, the x-axis is restricted to $x/n \leq 0.5$. (Recall that all HCs considered here are either restricted to $x/n \leq 0.5$ or symmetric around $x/n = 0.5$. Recall too that $n$ and $x$ are treated as continuous here.) For each E and each HC, the maximum value of the iso-E curve occurs at the TrP, with the segment to the left corresponding to evidence for $\theta_1$ and the segment to the right corresponding to evidence for $\theta_2$.

One way to use these graphs is to find the sample size $n$ corresponding to a particular value of E for given $x/n$. For instance, to obtain E = 2 in favor of $\theta_1$, we would need $n$ = 1.5, 1.1, 3.0 and 3.6 heads in a row ($x/n = 0$), for Class I(a), Class I(b), Class II(a) and Class II(b), respectively. Apparently E = 2 is quite easy to achieve, in the sense that relatively few tosses will yield E = 2 if they are all heads. By contrast, to get E = 4 one would need 7.0, 3.0, 15.2 and 20.5 tosses, all heads, for the four HCs respectively; while E = 8 (not shown in Figure) would require 21.8, 7.0, 67.3 and 106.6 heads, respectively. Another way to use the graphs is to see the "effect size" at which a given sample size $n$ will return evidence E. As Figure 6 shows, whether the evidence favors $\theta_1$ (left of TrP) or $\theta_2$ (right of TrP), much larger samples are required to achieve a given E the closer $x/n$ is to the TrP, or in other words, the less incompatible the data are with the non-favored hypothesis. For instance, for Class II(a) and considering evidence for $\theta_1$, for n = 100, E = 4 for $x/n \approx 0.07$; but for $n$ = 300, that same E = 4 is achieved for $x/n \approx 0.25$, a much smaller deviation from ½.

**Figure 6** Iso-E profiles comparing four HCs, for (a) E = 2, (b) E = 4. For purposes of illustration, $0.4 \leq \theta_2 \leq 0.6$ for Class II(b).

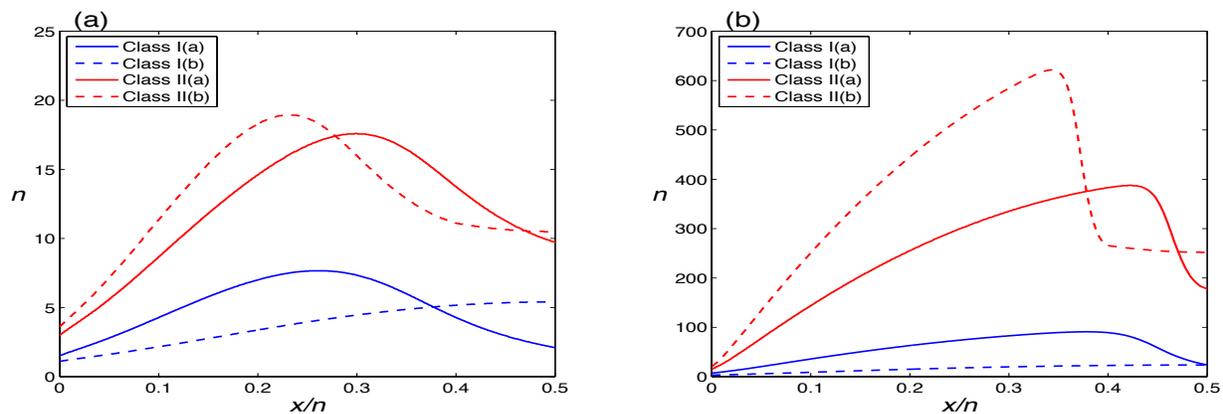

To our knowledge, ours is the only framework that generates a rigorous mathematical definition of what it means for evidence to be constant across different sets of data and different forms of HC. Note too that E is on a proper ratio scale [4, 17], so that 6(b) represents a doubling of the strength of evidence as shown in 6(a) (and E = 8 represents a doubling again relative to E = 4). This is a unique feature of E compared to all other proposed evidence measures of which we are aware. Figure 6 is a type of graph that can be meaningfully produced *only* once one has a properly calibrated measurement E in hand.

**Discussion**

With the results presented above, we have taken important steps forward towards generalizing our original information-dynamic theory in support of a properly calibrated measure E of statistical evidence. Three new results in particular move the theory forward. First, we have shown how to modify the original EqS (1.4) for one-sided HCs to obtain a new EqS (2.3) which handles two-sided HCs. More generally, we have shown that these two equations alone cover both non-nested and nested HCs, including a broad class of composite vs. simple or composite vs. composite comparisons. Second, while the original EqS had the same form as the ideal gas equation, the revised EqS needed to properly handle nested HCs has the same form as the thermodynamic Van der Waals equation. Third, we have discovered that the constant $c_1$, which corresponds to $c_V$ in the physical versions of these equations, seems to function in the information-dynamic equations as a kind of d.f. adjustment, allowing us for the first time to rigorously compare evidence across HCs of differing dimensionality.

To date we have focused on building this novel "plero"-dynamics (from the Greek word for information) methodology and understanding its relationship to thermo-dynamics. Because our motivation – proper measure-theoretic calibration of evidence – is distinct from the objectives of other schools of statistical thought, we have found it challenging to try to relate plerodynamics to components of standard mathematical statistical theory. But what emerges from the current results is a novel concept of evidence as a relationship between the maximum log LR and the area (or more generally, volume) under the LR, where the relationship is mediated by a quantity related to the Fisher information (for nested HCs, see Appendix 2), and also by something related to statistical degrees of freedom. This strongly suggests that we should be able to tie current results back to fundamental statistical theory. This will entail a detailed consideration of the concept of degrees of freedom, as it appears in plerodynamics, with its corresponding role in familiar statistical theory, and/or with its role in physical theory.

We had originally thought that every statistical model would require discovery of a separate EqS. But we now speculate that the EqS may depend *only* on the form of the HC, and not on the particular form of the likelihood. That is, our basic equations of state for the binomial model may extend to more complex models, based on general properties of likelihood ratios, at least under broad regularity conditions. Of course, so far the equations remain restricted to single-parameter models, and a somewhat restricted class of HCs (excluding "asymmetric" and non-partitioning HCs, as described above). We also have not considered extensions to continuous distributions. However, we follow Baskurt and Evans [18] in considering all applications of statistical inference as fundamentally about discrete, rather than continuous, distributions. This also raises the possibility of another way of relating plerodynamics back to thermodynamics, since in this case plerodynamics in its most general form could perhaps be represented solely in terms of the Boltzmann distribution.

In this paper we have not focused on the "-dynamic" part of plerodynamics, but the underlying theory motivating the approach taken here is very closely aligned with the macroscopic description of thermodynamic systems in terms of conservation and inter-conversion of heat and work. As we have noted previously, there is, however, no direct mapping of the basic thermodynamic variables (volume, pressure, mechanical work, number of particles) onto corresponding statistical variables. For example, the number of observations $n$ on the statistical side does *not* function in our EqS as the analogue of the number of particles in physics; rather, the number of observations $n$ appears to be part of the description of the statistical system's information "energy," rather than its size. (See [5] for discussion of this issue.) Therefore, but perhaps quite counter-intuitively, we do *not* expect to see a simple alignment of plerodynamics with statistical mechanical (microscopic) descriptions of physical systems, even in the event that we are able ultimately to consolidate the theory under the umbrella family of Boltzmann distributions.

It remains an open question how deep the connection between plero-and thermo-dynamics really runs. Our discovery here that a simple revision to the ideal gas equation of state solves one of the difficulties we have faced until now – our inability to generalize from a one-side to a two-sided hypothesis contrast – goes some distance towards vindicating our original co-opting of that particular equation of state for statistical purposes. The further discovery that the revised equation is identical in form to the Van der Waals equation surely takes us some distance further.

## Appendix 1: Maximum log LR plays the role of entropy, not evidence

We consider here the idea of the maximum log LR as an entropy term. We continue to restrict attention to the binomial likelihood in θ under the HCs considered in the main text. We note first that the maximum log LR is equivalent to a particular form of KL divergence (KLD), where

$$KLD[P_n(x;\,\theta_1), P_n(x;\,\theta_2)] \quad = \sum_{x=0}^{n} P_n(x;\,\theta_1) \log \frac{P_n(x;\,\theta_1)}{P_n(x;\,\theta_2)}$$

$$= \sum_{x=0}^{n} P_n(x;\,\theta_1) \log \frac{L(\theta_1; n, x)}{L(\theta_2; n, x)}$$

$$= \log \left(\frac{\theta_1}{\theta_2}\right) \sum_{x=0}^{n} x P_n(x;\,\theta_1) + \log \left(\frac{1-\theta_1}{1-\theta_2}\right) \sum_{x=0}^{n} (n-x) P_n(x;\,\theta_1)$$

$$= E[X] \log \left(\frac{\theta_1}{\theta_2}\right) + (n - E[X]) \log \left(\frac{1-\theta_1}{1-\theta_2}\right)$$

$$= n\theta_1 \log \left(\frac{\theta_1}{\theta_2}\right) + (n - n\theta_1) \log \left(\frac{1-\theta_1}{1-\theta_2}\right)$$

$$= n\theta_1 \log \theta_1 + (n - n\theta_1) \log(1-\theta_1) - n\theta_1 \log \theta_2 - (n - n\theta_1) \log(1 - \theta_2). \quad \text{(A1.1)}$$

If we now evaluate the KLD at $\theta_1 = \hat{\theta} = \frac{x}{n}$ and $\theta_2 = \hat{\theta}_i$, as in the main text, we have what we call the *observed* KL divergence ("observed" because the expectation is taken with respect to a probability distribution based on the data), which is equal to the log of the maximum likelihood ratio (MLR):

$$KLD_{OBS}[P_n(x;\,\hat{\theta}), P_n(x;\,\hat{\theta}_i)]$$

$$= x \log \left(\frac{x}{n}\right) + (n-x) \log \left(1 - \left(\frac{x}{n}\right)\right) - x \log \hat{\theta}_i - (n-x) \log(1 - \hat{\theta}_i)$$

$$= \max_{\theta} \log \frac{L(\theta; n, x)}{L(\hat{\theta}_i; n, x)} = \log \text{MLR}. \quad \text{(A1.2)}$$

Note that Kullback [12] and others [19] treat Kullback-Leibler divergence as a key quantity in an entropy-based inferential framework (see also [5]), while the MLR or its logarithm is sometimes interpreted as the statistical evidence for $\hat{\theta}$ against some simple alternative value [9, 20]. In our framework, the log MLR functions as the entropy term S (2.1).

There are several reasons why the MLR (or its logarithm) cannot be an evidence measure. First, as noted in the main text, the MLR violates important BBPs. In particular, when more maximization is done in the numerator than in the denominator, MLR $\geq 1$, and it cannot indicate evidence in favor of $\hat{\theta}_i$ or accumulate increasing evidence in favor of $\hat{\theta}_i$ as a function of increasing *n*, which violates elements of BBP(i)-(iii). Additionally, for fixed *x/n*, the MLR itself increases exponentially in *n*, while the log MLR increases linearly in *n*, both of which violate BBP(iv). (Indeed, the simple vs. simple LR itself, which is sometimes used to *define* the evidence [9, 20, 21], violates BBP(iv).) Yet there is clearly a reason why the MLR seems to function as a reasonably good proxy for an evidence measure under many circumstances.

We interpret the log MLR as the difference in *information* provided by the data for $\hat{\theta}$ vs. $\hat{\theta}_i$. As a general rule, as information (in an informal sense) goes up, so too does evidence. But information and evidence also must be distinguished, in the sense that increasing the amount of information might reduce the evidence for bias, if the more we toss the coin the closer $\hat{\theta}$ moves towards ½. Apparently evidence requires us to take account of information in the sense of KLD$_{\text{OBS}}$, or equivalently, in the form of the MLR, but not *only* information in this sense.

## Appendix 2: Calculation of *b*

Here we describe the rationale for setting the constant *b* as it appears in the main text. The BBPs impose severe constraints on the set of available solutions for *b*, and it appears that there is little leeway in choosing a functional form for *b* that allows us to express E for both Class II(a) and Class II(b) through a single EqS while maintaining the BBPs.

By experimentation (informed trial and error), we arrived at the following definition, which incorporates two rate constants: $r_1$, which controls the curvature of *b* over $\Theta_2$; and $r_2$, which controls the baseline value of *b* at the boundaries of this region, that is, at the points $\theta_{2l}$, $\theta_{2r}$. Let the value of *b* at these points be $b(\theta_{2l}) = b(\theta_{2r})$. We note up front that for given *n*, the minimum value of the Fisher information, Min FI$(n) = -E\left[\frac{d^2}{d\theta^2}\log L(\theta)\right]$, occurs when $\theta = $ ½. Then we have

$$b = \begin{cases} r_1V - r_2 \dfrac{\sqrt{2\pi}}{\sqrt{Min\,FI(n)}}; & x/n \in \Theta_2 \\ g\left[b(\theta_{2j}),0\right] & ;\,otherwise \end{cases} \tag{A2.1}$$

where $j = (l, r)$ and $g$ is the linear function connecting the points $b(\theta_{2l})$ and 0 (on the left) or $b(\theta_{2r})$ and 1 (on the right).

We set $r_1 = 2 - [\theta_{2r} - \theta_{2l}]$, so that the curvature of $b$ depends on the width of the $\Theta_2$ interval. We found that we needed to constrain $r_2$ such that $\tfrac{1}{2} \le (\,r_1 - (\tfrac{1}{2})r_2) \le \tfrac{3}{4}$. Thus we used $r_2 = 2r_1 - \tfrac{1}{2}(2 + [\theta_{2r} - \theta_{2l}])$. Figure 7 shows $b$ and V for various $\Theta_2$ for $n = 50$.

**Figure 7** Relationship among V, *b* and V-*b* using (A2.1) to calculate b. Shown here are four $\Theta_2$ intervals: (a) [0.49, 0.51], (b) [0.4, 0.6], (c) [0.3, 0.7], (d) [0.2, 0.8].

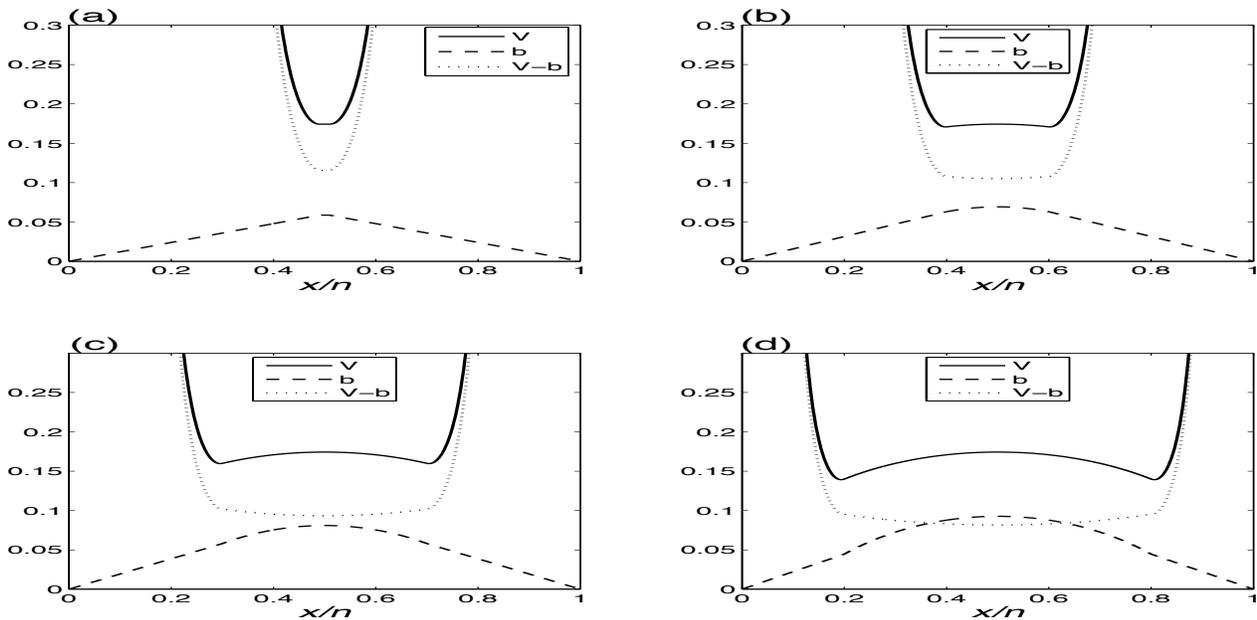

**Acknowledgments** This work was supported by a grant from the W.M. Keck Foundation. We thank Bill Stewart and Susan E. Hodge for helpful discussion, and we are indebted to Dr. Hodge for her careful reading of earlier drafts of this manuscript and for her many helpful comments, which led to substantial improvements in the paper.

**Author Contributions** Both authors contributed equally to this work.

**Conflict of Interests** The authors declare no conflict of interest.